\documentclass[reqno]{amsart}
\usepackage{hyperref}
\usepackage{amsmath}
\usepackage{amssymb}
\usepackage{amscd}
\usepackage{graphics}
\usepackage{latexsym}
\usepackage{stmaryrd}

\theoremstyle{plain}
\newtheorem{theorem}{Theorem}[section]
\newtheorem{thm}[theorem]{Theorem}
\newtheorem{cor}[theorem]{Corollary}
\newtheorem{lem}[theorem]{Lemma}
\newtheorem{prop}[theorem]{Proposition}

\theoremstyle{definition}

\newtheorem{ques}[theorem]{Question}

\theoremstyle{remark}

\newcommand{\marpar}[1]{}

\newcommand{\QQ}{\mathbb{Q}}

\newcommand{\CC}{\mathbb{C}}

\newcommand{\PP}{\mathbb{P}}

\newcommand{\mc}{\mathcal}

\newcommand{\OO}{\mc{O}}

\newcommand{\kgnb}[1]{\overline{\text{M}}_{#1}}

\newsavebox{\sembox}
\newlength{\semwidth}
\newlength{\boxwidth}

\newsavebox{\semrbox}
\newlength{\semrwidth}
\newlength{\boxrwidth}

\title[Linear sections]
{Fano varieties and linear sections of hypersurfaces}

\author[Starr]{Jason Michael Starr}
\address{Department of Mathematics \\
  Massachusetts Institute of Technology \\ Cambridge MA 02139}
\email{jstarr@math.mit.edu} 

\date{\today}
\begin{document}


\begin{abstract}
When $n$ satisfies an inequality which is almost best possible,
we prove that the $k$-plane sections of 
every smooth, degree $d$, complex hypersurface in $\PP^n$
dominate the moduli space of degree $d$ hypersurfaces in $\PP^k$.
As a corollary
we prove that, for $n$ sufficiently large, every smooth, degree $d$
hypersurface in $\PP^n$ satisfies a version of ``rational
simple connectedness''.  
\end{abstract}


\maketitle

\section{Statement of results}

In their article \cite{HMP}, Harris, Mazur and Pandharipande prove
that for fixed integers $d$ and $k$, there exists an integer
$n_0=n_0(d,k)$ such that for every $n\geq n_0$, 
every smooth degree $d$
hypersurface $X$ in $\PP^n_{\CC}$ has a number of good properties:
\begin{enumerate}
\item[(i)]
The hypersurface is unirational.
\item[(ii)]
The Fano variety of $k$-planes in $X$ has the expected dimension.
\item[(iii)]
The $k$-plane sections of the hypersurface dominate the moduli space
of degree $d$ hypersurfaces in $\PP^k$.
\end{enumerate}
It is this last property which we consider.
To be precise, the statement is that the following rational transformation
$$
\Phi: \mathbb{G}(k,n) \dashrightarrow \PP^{N_d}//\mathbf{PGL}_{k+1}
$$
is dominant.  Here $\mathbb{G}(k,n)$ is the Grassmannian parametrizing
linear $\PP^k$s in $\PP^n$, $\PP^{N_d}$ is
the parameter space for degree $d$ hypersurface in $\PP^k$,
$\PP^{N_d}//\mathbf{PGL}_{k+1}$ is the moduli space of semistable
degree $k$ hypersurface in $\PP^k$, and $\Phi$ is the rational
transformation sending a $k$-plane $\Lambda$ to the moduli point of
the hypersurface $\Lambda \cap X \subset \Lambda$ (assuming
$\Lambda\cap X$ is a semistable degree $k$ hypersurface in $\PP^k$).  

The bound $n_0(d,k)$ is very large, roughly a $d$-fold iterated
exponential.  Our result is the following.

\begin{thm} \label{thm-main}
\marpar{thm-main}
Let $X$ be a smooth degree $d$ hypersurface in $\PP^n$.  The map
$\Phi$ is dominant if 
$$
n\geq \binom{d+k-1}{k} + k - 1.
$$
\end{thm}

\begin{ques} \label{ques-n}
\marpar{ques-n}
For fixed $d$ and $k$, what is the
smallest integer $n_0 = n_0(d,k)$ such that for every $n\geq n_0$ and
every smooth, degree $d$ hypersurface in $\PP^n$, the associated
rational transformation $\Phi$ is dominant?  
\end{ques}

Theorem ~\ref{thm-main} is equvialent to the inequality
$$
n_0(d,k) \leq \binom{d+k-1}{k} + k - 1.
$$
If $\Phi$ is dominant, then the dimension of the domain is at least
the dimension of the target, i.e., 
$$
(k+1)(n-k) = \text{dim}\mathbb{G}(k,n) \geq
\text{dim}(\PP^{N_d}//\mathbf{PGL}_{k+1}) = \binom{d+k}{k} - (k+1)^2.
$$
This is equivalent to the condition
$$
n_0(d,k) \geq \frac{1}{k+1}\binom{d+k}{k} - 1.
$$
As far as we know, this is the correct bound.  The bound from
Theorem~\ref{thm-main} differs from this optimal bound by roughly a
factor of $k$.

The main step in the proof is a result of some independent interest.

\begin{prop} \label{prop-Fano}
\marpar{prop-Fano}
Let $X$ be a smooth degree $d$ hypersurface in $\PP^n$.  Let $F_k(X)$
be the Fano variety of $k$-planes in $X$.  There exists an irreducible
component $I$ of $F_k(X)$ of the expected dimension if
$$
n\geq \binom{d+k-1}{k} + k.
$$
Moreover, if
$$
n = \binom{d+k-1}{k} + k - 1
$$
then there is a nonempty open subset $U_{k-1} \subset F_{k-1}(X)$ such
that for every $[\Lambda_{k-1}]\in U_{k-1}$, there exists no
$k$-plane in $X$ containing $\Lambda_{k-1}$.
\end{prop}

Theorem ~\ref{thm-main} implies a result about rational curves on
every smooth hypersurface of sufficiently small degree.  The
Kontsevich moduli space $\kgnb{0,r}(X,e)$ parametrizes isomorphism
classes of data
$(C,q_1,\dots,q_r,f)$ of a proper, connected, at-worst-nodal,
arithmetic genus $0$ curve $C$, an ordered collection $q_1,\dots,q_r$
of distinct smooth points of $C$ and a morphism $f:C\rightarrow X$
satisfying a stability condition.  The space $\kgnb{0,r}(X,e)$ is
projective.  There is an evaluation map
$$
\text{ev}:\kgnb{0,r}(X,e) \rightarrow X^r
$$
sending a datum $(C,q_1,\dots,q_r,f)$ to the ordered collection
$(f(q_1),\dots,f(q_r))$.  

\begin{cor} \label{cor-main}
\marpar{cor-main}
Let $X$ be a smooth degree $d$ hypersurface in $\PP^n$.  If
$$
n \geq \binom{d^2+d-1}{d-1} + d^2 -1
$$
then for every integer $e\geq 2$ there exists a canonically defined
irreducible component $\mc{M} \subset \kgnb{0,2}(X,e)$ such that the
evaluation morphism
$$
\text{ev}: \mc{M} \rightarrow X \times X
$$
is dominant with rationally connected generic fiber, i.e., $X$
satisfies a version of
rational simple connectedness.  Moreover $X$ has a very twisting family
of pointed lines, cf. ~\cite[Def. 3.7]{Sr1c}.
\end{cor}

This is proved in \cite{Sr1c} assuming $n$ satisfies a
much weaker hypothesis
$$
n \geq d^2
$$
but only for \emph{general} hypersurfaces, not for \emph{every} smooth
hypersurface.  
The goal here is to find a stronger hypothesis on $n$ that
guarantees the theorem for every smooth hypersurface.  

\section{Flag Fano varieties} \label{sec-FF}
\marpar{sec-FF}

Naturally enough, the proof of Proposition~\ref{prop-Fano} uses an
induction on $k$.  To set up the induction it is useful to consider
not just $k$-planes in $X$, but flags of linear spaces 
$$
\PP^0 \subset \PP^1 \subset \PP^2 \subset \dots \subset \PP^k \subset
X.
$$
The variety parametrizing such flags is the \emph{flag Fano variety}
of $X$.  
Also, although
we are ultimately interested only in the case of a hypersurface
in projective space, for the induction it is useful to allow a
more general projective subvariety.

Let $S$ be a scheme such that $H^0(S,\OO_S)$ contains $\QQ$.  
Let $E$ be a locally free $\OO_S$-module of rank
$n+1$, and let $X \subset \PP E$ be a 
closed subscheme such that the
projection $\pi:X\rightarrow S$ is smooth and surjective
of constant relative
dimension $\text{dim}(X/S)$.  
In other words, $X$ is a
family of smooth, $\text{dim}(X/S)$-dimensional
subvarieties of $\PP^n$ parametrized by $S$.  
  
Let $0\leq k \leq n$ be an
integer.  Denote by $\text{Fl}_k(E)$ the partial flag manifold representing
the functor on $S$-schemes
$$
T \mapsto \{(E_1 \subset E_2 \subset \dots \subset E_{k+1} \subset
E_T)| E_i \text{ locally free of rank } i, i=1,\dots,k+1\}.
$$
For every $0 \leq j \leq k \leq n$, denote by $\rho_k^j:\text{Fl}_k(E)
\rightarrow \text{Fl}_j(E)$ the obvious projection.  The \emph{flag
  Fano variety} is the locally closed subscheme
$\text{Fl}_k(X) \subset \text{Fl}_k(E)$ 
parametrizing flags such that $\PP(E_{k+1})$ is contained in $X$.  In
particular, $\text{Fl}_0(X) = X$.  Denote by $\rho_k^j:\text{Fl}_k(X)
\rightarrow \text{Fl}_j(X)$ the restriction of $\rho_k^j$.

\subsection{Smoothness} \label{subsec-sm}
\marpar{subsec-sm}

There are two elementary observations about the schemes $\text{Fl}_k(X)$.

\begin{lem} \label{lem-obs1}\cite[1.1]{KMM92c}
\marpar{lem-obs1}
There exists an open dense subset $U\subset X$ such that $U\times_X
\text{Fl}_1(X)$ is smooth over $U$.
\end{lem}

\begin{lem} \label{lem-obs2} \marpar{lem-obs2} 
  Set $S^\text{new} = U$, the open subset from Lemma~\ref{lem-obs1}.
  Set $E^\text{new}$ to be the universal rank $n$ quotient bundle of
  $\pi^* E|_U$ so that $\PP(E^\text{new}) = U\times_{\PP(E)}
  \text{Fl}_1(E)$ and set $X^\text{new} = \text{Fl}_1(U)$.  Then for
  every $0\leq k \leq n-1$, $\text{Fl}_k(X^\text{new}) = U\times_X
  \text{Fl}_{k+1}(X)$.
\end{lem}

\begin{proof}
This is obvious.
\end{proof}

\begin{prop} \label{prop-U}
\marpar{prop-U}
There exists a sequence of open subschemes $(U_k \subset
\text{Fl}_k(X))_{0\leq k \leq n}$ satisfying the following conditions.
\begin{enumerate}
\item[(i)]
The open subset $U_0$ is dense in $\text{Fl}_0(X)$, and for every $1\leq k
\leq n$, $U_k$ is 
dense in $(\rho_k^{k-1})^{-1}(U_{k-1})$. 
\item[(ii)]
For every $1\leq k \leq n$, $\rho_k^{k-1}:
(\rho_k^{k-1})^{-1}(U_{k-1})
\rightarrow U_{k-1}$ is
smooth.  
\end{enumerate}
\end{prop}

\begin{proof}
  Let $U_0$ be the open subscheme from Lemma~\ref{lem-obs1}.  By way
  of induction, assume $k>0$ and the open subscheme $U_{k-1}$ has been
  constructed.  As in Lemma~\ref{lem-obs2}, replace $S$ by $U_{k-1}$,
  replace $E$ by the universal quotient bundle, and replace $X$ by
  $(\rho_k^{k-1})^{-1}(U_{k-1})$.  Now define $U_k \subset
  (\rho_k^{k-1})^{-1}(U_{k-1})$ to be the open subscheme from
  Lemma~\ref{lem-obs1}.
\end{proof}

\subsection{Dimension} \label{subsec-dim}
\marpar{subsec-dim}

Using the Grothendieck-Riemann-Roch formula, it is possible to express
the Chern classes of $U\times_X F_1(X)$ in terms of the Chern classes of
$U$.  Iterating this leads, in particular, to a formula for the
dimension of $U_k$.  Denote by $G_1$, resp. $G_2$, the restriction to
$\text{Fl}_1(U)$ of $E_1$, resp. $E_2$.  Denote by $L$ the invertible
sheaf
$$
L := (G_2/G_1)^\vee.
$$
Denote by
$$
\pi:\PP G_2 \rightarrow \text{Fl}_1(U),
$$
$$
\sigma: \text{Fl}_1(U) = \PP G_1 \rightarrow \PP G_2,
$$
and
$$
f: \PP G_2 \rightarrow X
$$
the obvious morphisms.  In other words, $\PP G_2$ is a family of
$\PP^1$s over $\text{Fl}_1(U)$, $\sigma$ is a marked point on each
$\PP^1$, and $f$ is an embedding of each $\PP^1$ as a line in $X$.
The formula for the Chern character of the vertical tangent bundle of
$\rho_1^0$ is,
$$
\text{ch}(T_{\text{Fl}_1(U)/U}) = \pi_*f^*[(\text{ch}(T_{X/S}) -
\text{dim}(X/S))\text{Todd}(\OO_{\PP E}(1)|_X)] - \text{ch}(L) - 1.
$$
Given a flag $\PP = (\PP^1\subset \PP^2 \subset \dots \subset \PP^k \subset
\PP^n)$ in $U_k$, the formula for the fiber dimension of
$\rho_k^{k-1}$ at $\PP$ is 
$$
\text{dim}(U_k/U_{k-1}) = \sum_{m=1}^k b_{k,m}\langle
\text{ch}_m(T_{X/S}),\PP^m\rangle -k-1
$$
where $\text{ch}_m(E)$ is the $m^\text{th}$ graded piece of the Chern
character of $E$, and where the coefficients $b_{k,m}$ are the unique
rational numbers such that
$$
\binom{x+k-1}{k} = \sum_{m=1}^k \frac{b_{k,m}}{m!} x^m.
$$

Now define the numbers $a_{k,m}$ to be
$$
a_{k,m} = \sum_{l=m}^k b_{l,m},
$$
in other words,
$$
\sum_{m=1}^k \frac{a_{k,m}}{m!}x^m = \sum_{l=1}^k \binom{x+l-1}{l}.
$$
Then it follows from the previous formula that the dimension of $U_k$
at $\PP$ equals
$$
\text{dim}(U_k) = \sum_{m=1}^k a_{k,m}\langle
\text{ch}_m(T_{X/S}),\PP^m\rangle + \text{dim}(X) - k^2.
$$

In a related direction, there is a class of complex projective
varieties that
is stable under the operation of replacing $X$ by a general fiber of
$\text{Fl}_1(X) \rightarrow X$.  Call a subvariety $X$ of $\PP^n$ a
\emph{quasi-complete-intersection} of type
$$
\underline{d} = (d_1,\dots,d_c)
$$ 
if there is a sequence
$$
X = X_c \subset X_{c-1} \subset \dots
\subset X_1 \subset X_0 = \PP^n
$$ 
such that each $X_k$ is a Cartier divisor in $X_{k-1}$ in the linear
equivalence class of $\OO_{\PP^n}(d_k)|_{X_{k-1}}$.  If $X$ is a
quasi-complete-intersection, then every fiber of $U\times_X
\text{Fl}_1(X) \rightarrow U$ is also a quasi-complete-intersection in
$\PP^{n-1}$ of type
$$
(1,2,\dots,d_1,1,2,\dots,d_2,\dots,1,2,\dots,d_c).
$$  
Iterating this, every (non-empty) fiber of
$(\rho_k^{k-1})^{-1}(U_{k-1}) \rightarrow U_{k-1}$ is a
quasi-complete-intersection in $\PP^{n-k}$ of dimension
$$
N_k(n,\underline{d}) = n - k - \sum_{i=1}^c \binom{d_i+k-1}{k}.
$$
Since the $m^\text{th}$ graded piece of the Chern character of $T_X$
equals
$$
\text{ch}_m(T_X) = (n+1 - \sum_{i=1}^c d_i^m)c_1(\OO(1))^m / m!
$$
this agrees with the previous formula for the fiber dimension.
 
\begin{cor} \label{cor-sm}
\marpar{cor-sm}
Let $X$ be a smooth quasi-complete-intersection of type
$\underline{d}$.  If the integer $N_k(n,\underline{d})$ is
nonnegative, there exists an irreducible component $I$ of
$\text{Fl}_k(X)$ having the expected dimension
$$
\text{dim}(I) = \sum_{m=0}^k N_m(n,\underline{d}).
$$
\end{cor}

\begin{proof}
  Of course we define $I$ to be the closure of any connected component
  of $U_k$.  The issue is whether or not $U_k$ is empty.  By
  construction $U_k$ is not empty if for every $m=1,\dots,k$ the
  morphism $\rho_m^{m-1}$ is surjective.  By the argument above every
  fiber of $\rho_m^{m-1}$ is an iterated intersection in $\PP^{n-m}$
  of pseudo-divisors (in the sense of \cite[Def. 2.2.1]{F}) in the
  linear equivalence class of an ample divisor.  Thus the fiber is
  nonempty if the number of pseudo-divisors is $\leq n-m$.  This
  follows from the hypothesis that $N_k(n,\underline{d}) \geq 0$.
\end{proof}

\section{Proofs} \label{sec-proofs}
\marpar{sec-proofs}

\begin{proof}[Proof of Proposition ~\ref{prop-Fano}]
The first part follows from Corollary ~\ref{cor-sm}.  For the second
part, observe that if $N_k(n,d) = -1$, then $N_{k-1}(n,d)$ is
nonnegative.  Therefore, by the first part, the open subset
$U_{k-1}$ from Proposition \ref{prop-U} is nonempty.
Since $(\rho_k^{k-1})^{-1}(U_{k-1}) \rightarrow U_{k-1}$ is smooth of
the expected dimension, and since the expected dimension is negative,
$(\rho_k^{k-1})^{-1}(U_{k-1})$ is empty.  In other words, for every
$[\Lambda_{k-1}] \in U_{k-1}$, there exists no $k$-plane in $X$ containing
$\Lambda_{k-1}$.  
\end{proof}

\begin{proof}[Proof of Theorem~\ref{thm-main}]
Let $(H_{k,n},e)$ be the universal pair of a scheme $H_{k,n}$ and a
closed immersion of $H_{k,n}$-schemes
$$
(\text{pr}_H,e):H_{k,n} \times \PP^k \rightarrow H_{k,n} \times \PP^n
$$
whose restriction to each fiber $\{h\} \times \PP^k$ is a linear
embedding.  In other words, $H_{k,n}$ is the open subset of $\PP
\text{Hom}(\CC^{k+1},\CC^{n+1})$ parametrizing injective matrices.  Of
course there is a natural action of $\mathbf{PGL}_{k+1}$ on $H_{k,n}$,
and the quotient is the Grassmannian $\mathbb{G}(k,n)$.  Denote by
$\widetilde{F}_k(X)$ the inverse image of $F_k(X)$ in $H_{k,n}$, i.e.,
$\widetilde{F}_k(X)$ parametrizes linear embeddings of $\PP^k$ into
$X$.  

Let $F$ be a defining equation for the hypersurface $X$.  Then $e^* F$
is a global section of $e^*\OO_{\PP^n}(d)$.  By definition, this is
canonically isomorphic to $\text{pr}_{\PP^k}^* \OO_{\PP^k}(d)$.
Therefore $e^* F$ determines a regular morphism
$$
\widetilde{\Phi}:H_{k,n} \rightarrow H^0(\PP^k,\OO_{\PP^k}(d)).
$$
Denote by $V$ the open subset of $H_{k,n}$ of points whose fiber
dimension equals 
$$
\text{dim}H_{k,n} - \text{dim} H^0(\PP^k,\OO_{\PP^k}(d)).
$$
The rational transformation $\Phi$ is dominant if and only if
$\widetilde{\Phi}$ is dominant.  And 
the morphism $\widetilde{\Phi}$ is dominant if and only if $V$ is
nonempty.  

The scheme $\widetilde{F}_k(X)$ is the fiber
$\widetilde{\Phi}^{-1}(0)$.  
If
$$
n\geq \binom{d+k-1}{k} + k
$$
then Proposition ~\ref{prop-Fano} implies there exists an irreducible
component $I$ of $F_k(X)$ of the expected dimension.  Thus the inverse
image $\widetilde{I}$ in $H_{k,n}$ is an irreducible component of
$\widetilde{F}_k(X)$ of the expected dimension, or what is equivalent,
the expected codimension.  But the expected codimension is precisely
$$
h^0(\PP^k,\OO_{\PP^k}(d)) = \binom{d+k}{k}.
$$
Thus, the generic point of $\widetilde{I}$ is contained in $V$, i.e.,
$V$ is not empty.

This only leaves the case when
$$
n= \binom{d+k-1}{k} + k-1.
$$
The argument is very similar.  Let $y$ be a linear coordinate on
$\PP^k$, and let $\widetilde{G}_k(X)$ be the closed subscheme of
$H_{k,d}$ where $e^* F$ is a multiple of $y^d$.  In other words,
$\widetilde{G}_k(X)$ parametrizes linear embeddings of $\PP^k$ into
$\PP^n$ whose intersection with $X$ contains $d\mathbb{V}(y)$. 
There
is a projection morphism $\widetilde{G}_k(X) \rightarrow F_{k-1}(X)$
associating to the linear embedding the $(k-1)$-plane 
$$
\Lambda_{k-1} = \text{Image}(\mathbb{V}(y)).
$$
Denote by $G_k(X)$ the image of $\widetilde{G}_k(X)$ under the obvious
morphism
$$
\widetilde{G}_k(X) \rightarrow F_{k-1}(\PP^n)\times F_k(\PP^n).
$$

Recall that for a quasi-complete-intersection $X$, the fiber of
$F_1(X) \rightarrow X$ is an interated intersection of ample
pseudo-divisors in projective space.
By a very similar argument, every fiber of
$G_k(X) \rightarrow F_{k-1}(X)$ is an iterated intersection
of ample pseudo-divisors in the
projective space $\PP^n/\Lambda_{k-1} \cong \PP^{n-k}$.  Moreover, the
fiber of $\text{Fl}_k(X) \rightarrow \text{Fl}_{k-1}(X)$ (for any
extension of $\Lambda_{k-1}$ to a flag in $\text{Fl}_{k-1}(X)$)
is an ample pseudo-divisor in $G_k(X)$.  
By the second part of Proposition \ref{prop-Fano}, there exists a
nonempty open subset $U_{k-1} \subset \Lambda_{k-1}$ such that for every
$\Lambda_{k-1} \in U_{k-1}$ this ample pseudo-divisor is empty.
Therefore the fiber in $G_k(X)$ is finite or empty.  But the equation
$$
n = \binom{d+k-1}{k} + k - 1
$$
implies the expected dimension of the fiber is $0$.  Since an
intersection of ample pseudo-divisors is nonempty if the expected
dimension is nonnegative,
the
fiber of $G_k(X) \rightarrow F_{k-1}(X)$ 
is not empty and has the expected dimension $0$.  Since $U_{k-1}$ has
the expected dimension, the open set $U_{k-1}\times_{F_{k-1}(X)}
\widetilde{G}_k(X)$ is nonempty and has the expected dimension.  Thus
it has the expected codimension.  Therefore a generic point of this
nonempty open set is in $V$, i.e., $V$ is not empty.
\end{proof}

\begin{proof}[Proof of Corollary~\ref{cor-main}]
Let $\mc{M}_e$ be an
irreducible component of $\kgnb{0,0}(X,e)$ not entirely contained in
the boundary $\Delta$.  Then for every integer $r\geq 0$ there exists
a unique irreducible component $\mc{M}_{e,r}$ of $\kgnb{0,r}(X,e)$
whose image in $\kgnb{0,0}(X,e)$ equals $\mc{M}_e$.  
Before defining the irreducible component $\mc{M}$ of
$\kgnb{0,2}(X,e)$, we will first inductively define
an irreducible
component 
$\mc{M}_e$ of $\kgnb{0,0}(X,e)$ which is not entirely contained in
the boundary $\Delta$ and such that the evaluation morphism
$$
\text{ev}:\mc{M}_{e,1} \rightarrow X
$$ 
is surjective.
Then we define $\mc{M}$ to be $\mc{M}_{e,2}$.

Let $U$ denote the
open subset of $\kgnb{0,1}(X,1)$ where the evaluation morphism
$$
\text{ev}:\kgnb{0,1}(X,1) \rightarrow X
$$
is smooth, i.e., $U$ parametrizes \emph{free} pointed lines.  By
\cite[1.1]{KMM92c}, $U$ contains every general fiber of $\text{ev}$.
By the argument in Subsection ~\ref{subsec-dim} (or any number of
other references), a general fiber of $\text{ev}$ is connected if $d
\leq n-2$.  Therefore $U\times_X U$ is irreducible.  There is an
obvious morphism $U\times_X U \rightarrow \kgnb{0,0}(X,2)$.  By
elementary deformation theory, the morphism is unramified and
$\kgnb{0,0}(X,2)$ is smooth at every point of the image.  Therefore
there is a unique irreducible component $\mc{M}_2$ of
$\kgnb{0,0}(X,2)$ containing the image of $U\times_X U$.  Because
$U\rightarrow X$ is dominant, $\mc{M}_2 \rightarrow X$ is also dominant.  

By way of induction assume $e\geq 3$ and $\mc{M}_{e-1}$ is given.
Form the fiber product $\mc{M}_{e-1,1}\times_X U$.  As above this is
irreducible, and there is an unramified morphism
$$
\mc{M}_{e-1,1}\times_X U \rightarrow \kgnb{0,0}(X,e)
$$
whose image is in the smooth locus.  Therefore there exists a unique
irreducible component $\mc{M}_e$ of $\kgnb{0,0}(X,e)$ containing the
image of $\mc{M}_{e-1,1}\times_X U$.  Because
$\mc{M}_{e-1,1}\rightarrow X$ is dominant, 
$\mc{M}_{e,1}\rightarrow X$ is also dominant.  This finishes the
inductive construction of the irreducible components $\mc{M}_e$, and
thus also of $\mc{M}_{e,2}$.  

It remains to prove that
$$
\text{ev}:\mc{M}_{e,2} \rightarrow X\times X
$$
is dominant with rationally connected generic fiber.  The article 
\cite{Sr1c}
gives an inductive argument for proving this.
To carry out the induction, one needs two results: 
the base of the induction and an important component of the induction
argument.  Set $k$ to be $d^2$.  For a general degree $d$ hypersurface
$Y$ in $\PP^k$, \cite[Prop. 4.6,
Prop. 10.1]{Sr1c} prove the two results for $Y$.  
By Theorem ~\ref{thm-main}, since
$$
n \geq \binom{d+k-1}{k} + k - 1,
$$
for a general $\PP^k \subset \PP^n$ the intersection $Y = \PP^k \cap
X$ is a general degree $d$ hypersurface in $\PP^k$.  Thus the two
results hold for $Y$.  As is clear from the proofs of
\cite[Prop. 4.6, Prop. 10.1]{Sr1c}, the results for $Y$ imply the
corresponding 
results for $X$.  
\end{proof}

\bibliography{my}

\begin{thebibliography}{1}

\bibitem{F}
W.~Fulton.
\newblock {\em Intersection theory}, volume~2 of {\em Ergebnisse der Mathematik
  und ihrer Grenzgebiete (3) [Results in Mathematics and Related Areas (3)]}.
\newblock Springer-Verlag, Berlin, 1984.

\bibitem{HMP}
J.~Harris, B.~Mazur, and R.~Pandharipande.
\newblock Hypersurfaces of low degree.
\newblock {\em Duke Math. J.}, 95(1):125--160, 1998.

\bibitem{KMM92c}
J.~Koll{\'a}r, Y.~Miyaoka, and S.~Mori.
\newblock Rational connectedness and boundedness of {F}ano manifolds.
\newblock {\em J. Differential Geom.}, 36(3):765--779, 1992.

\bibitem{Sr1c}
J.~Starr.
\newblock Hypersurfaces of low degree are rationally simply-connected.
\newblock preprint, 2004.

\end{thebibliography}
\bibliographystyle{abbrv}

\end{document}